\documentclass{amsart}

\usepackage{amscd,amssymb}

\usepackage{graphicx}
\usepackage{hyperref}

\usepackage{kpfonts}  
\usepackage[T1]{fontenc}
\usepackage[cp1250]{inputenc}
\usepackage[polish,english]{babel}

\newtheorem{lem}{Lemma}[section]
\newtheorem{cor}[lem]{Corollary}

\newtheorem{thm}[lem]{Theorem}
\newtheorem{prop}[lem]{Proposition}
\newtheorem{prty}[lem]{Property}
\newtheorem{fact}[lem]{Fact}

\makeatletter
\def\th@plain{%
  \thm@notefont{}
  \itshape 
}
\def\th@definition{%
  \thm@notefont{}
  \normalfont 
}
\def\th@procedure{
	\thm@notefont{}
	\ttfamily
}
\makeatother

\theoremstyle{definition} 
\newtheorem{ex}[lem]{Example}
\newtheorem{rk}[lem]{Remark}

\theoremstyle{procedure} 
\newtheorem{proc}{Procedure}

\newcommand{\Za}{\mathbb{Z}}

\newcommand{\re}{\mathbb{R}}
\newcommand{\ce}{\mathbb{C}}

\newcommand{\tr}{\triangle}
\newcommand{\supp}{{\rm supp}}

\author{Maria Michalska}
\author{Justyna Walewska}
\title[Milnor numbers of deformations]{Milnor numbers of deformations of semi-quasi-homogeneous plane curve singularities}
\address{Wydzia\l{} Matematyki i Informatyki, Uniwersytet \L{}\'o{}dzki, Banacha 22, 90-238 \L{}\'o{}d\'z{}, Poland}
\email{Maria.Michalska@math.uni.lodz.pl}
\email{walewska@math.uni.lodz.pl}
\keywords{Milnor numbers, deformations of singularities, nondegenerate singularities, Euclid's algorithm, Newton polygon}
\subjclass[2010]{14B07, 14N10, 32S30}

\begin{document}
\maketitle

\begin{abstract}
The aim of this paper is to show the possible Milnor numbers of deformations of semi-quasi-homogeneous isolated plane curve singularities. Main result states that if $f$ is irreducible and nondegenerate, by deforming $f$ one can attain all Milnor numbers ranging from $\mu(f)$ to $\mu(f)-r(p-r)$, where $r$ and $p$ are easily computed from the Newton diagram of $f$.

\end{abstract}

\section*{Introduction}

The main goal of this paper is to identify all possible Milnor numbers attained by deformations of plane curve singularities. This question is closely related to some of Arnold's problems~\cite{ArnoldProblems}. A~direct motivation for our study was a talk of Arkadiusz Płoski on recent developments and open questions regarding jumps of Milnor numbers given at the Łódź-Kielce seminar in June 2013 as well as questions from an article of Arnaud Bodin~\cite{Bo}. 

The most interesting point is establishing the initial Milnor jumps i.e. the greatest Milnor numbers attained by deformations. As was shown in general in~\cite{GuseinZade} and for special cases in~\cite{BK}, it is possible that not all Milnor numbers are attained, meaning that the jumps may be greater than one. Moreover, in these cases the Milnor numbers that are not attained give exactly the first jump greater than one. These results are related to bounds on Milnor numbers of singularities and refer to questions on possible Milnor numbers of singularities of given degree, see for instance \cite{Pl14} or~\cite{GreuelShustinL}. Moreover, the fact that the first jump is not equal one has in turn interesting implications for multiparameter versal deformations and adjacency of $\mu$-constant strata~\cite{ArnoldProblems}. In this paper we show, in particular, that this is not the case when considering irreducible nondegenerate semi-quasi-homogeneous singularities.

The approach presented here stems from the observation that many properties of the sequence of Milnor numbers attained by deformations of a singularity are possible to be established combinatorially, a fact that was not in our opinion sufficiently explored. In this paper we focus on the irreducible case which in our opinion is the hardest. A careful analysis shows that for semi-quasi-homogeneous singularities, assuming nondegeneracy, the problem boils down to three cases depending on the greatest common divisor of $p$ and $q$ (using the notation \eqref{eqPostacfQSH}). The irreducible case in such a setting is equivalent to saying that $p$ and $q$ are coprime. We show that initial jumps of Milnor numbers are equal to one. This result, on its own, can be used iteratively for many singularities to prove that all jumps are equal to one, as shown in Section~\ref{sectionRks}. On the other hand, we think of this paper as introduction to more general results based on the observation that if the procedure presented here is adjusted, it implies also solutions in general in the other two cases. For instance, given an isolated singularity $f$ of the form~\eqref{eqPostacfQSH} with ${\rm GCD}(p,q)=g$ such that $1<g<p$, one can show that the first jump is not bigger than $g$ (as was already shown in \cite{Bo} and \cite{W1}) but all Milnor numbers ranging from $\mu(f)-g$ to $\mu(f)- m-r(p-r)+1$ can be attained by deformations of $f$ (under notation $q\equiv r ({\rm mod}\ p)$). We defer the details to a subsequent publication. One would also like to note that parallel and complimentary research of the problem of jumps of Milnor numbers is a recent paper \cite{BKW}.

This article is organised as follows. First, we state the main result. In Section \ref{sectionPreliminaries} we begin with introducing notation that we hope will provide more clarity to further considerations. In paragraphs \ref{subsecNewtdiag} and \ref{subsecNewtonNumber} we recall some properties of the Newton diagram and Newton numbers. General combinatorial remarks and a reminder on Euclid's algorithm follow in paragraphs \ref{subsecGeneralComb}, \ref{subsecEEA} and \ref{subsecEEArks}. 

Section~\ref{secLemmas} presents steps needed in the proof of Theorem~\ref{mainThm}. It is divided into three parts. In Section~\ref{subsecDecrease} we prove validity of Procedure~\ref{procka1} that gives minimal jumps and allows to substitute $p$ and $q$ by $n(a-a')+a'$ and $n(b-b')+b'$ respectively (compare table \eqref{EEA}). In Section~\ref{subsecReductionLine} we prove an iteration of this procedure, that is Procedure~\ref{procka2}, is valid and gives minimal jumps until $p,q$ are reduced to $n'a'+a''$ and $n'b'+b''$ respectively. Whereas in Section~\ref{subsecShortEEA} we deal with the case (or the last line in Euclid's Algorithm) when $q\equiv \pm 1({\rm mod\ } p)$. 

Section~\ref{secMainThmComb} brings the procedures together to prove Theorem~\ref{thmMainCombinatorially}. The main Theorem~\ref{mainThm} follows immediately. The article concludes with some remarks and observations on further developments.

\section{Statement of the main result}

Throughout this paper we will consider an isolated plane curve singularity $f$ i.e.~the germ $f:(\ce^2,0)\to(\ce,0)$ is analytic and $0$ is the only solution of the system of equations $\nabla f(x,y)=f(x,y)=0$. By a deformation of $f$ we mean any analytic function $F:(\ce^{3},0)\to\ce$ such that $F(0,\cdot)=f$ and $F(t,\cdot)$ is an isolated singularity for every $t$ small enough.

The Milnor number $\mu(f)$ of an isolated singularity $f$ is the multiplicity of $\nabla f$ at zero. A classic result is that the Milnor number of a deformation $F$ of $f$ always satisfies the inequality $\mu(f)\ge \mu(F(t,\cdot))$ for $t$ small enough, see for instance \cite{GreuelShustinL}. Hence it makes sense to consider the strictly decreasing sequence $(\mu_i)_{i=0,\dots,w}$ of all positive integers attained as Milnor numbers of deformations of $f$. We have $\mu_0=\mu(f)$ and $\mu_w=1$. The sequence of positive integers $(\mu_{i-1}-\mu_{i})_{i=1,\dots,w}$ will be henceforth called the sequence of jumps of Milnor numbers.

We will consider the isolated singularity $f$ of the form
\begin{equation}
\label{eqPostacfQSH}
f=\sum_{p\alpha + q\beta\  \geq\  pq}c_{\alpha\beta}\ x^{\alpha}y^{\beta}
\end{equation} 
for some positive integers $p,q$.

\begin{thm}
\label{mainThm}
Given a nondegenerate isolated singularity $f$ of the form~\eqref{eqPostacfQSH} with $p<q$ coprime the sequence of Milnor jumps begins with
$$\underbrace{1\ ,\ \dots\ ,\ 1}_{r(p-r)}$$
where $r$ is the rest out of division of $q$ by $p$.
\end{thm}
\noindent{\bf Proof.} The proof follows immediately from Kouchnirenko's theorem (see Fact \ref{factNewMil}) and minimality of the jumps in Theorem \ref{thmMainCombinatorially}.\hfill$\blacksquare$

If $f$ is nondegenerate of the form~\eqref{eqPostacfQSH}, $p,q$ are coprime and $c_{p,0}c_{0,q}\neq 0$, then $f$ is irreducible. On the other hand, for any nondegenerate irreducible isolated singularity $f$, it is of the form \eqref{eqPostacfQSH}, $p,q$ are coprime and $c_{p,0}c_{0,q}\neq 0$. Hence as a special case a direct generalisation of \cite[Theorem 2]{Bo} follows. Namely

\begin{cor}
Given an irreducible nondegenerate isolated singularity, the claim of Theorem \ref{mainThm} holds.
\end{cor}

\section{Preliminaries on combinatorial aspects}
\label{sectionPreliminaries}

\subsection{Notations}

A Newton diagram of a set of points $\mathcal{S}$ is the convex hull of the set
$$\bigcup_{P\in\mathcal{S}} \left(P+\re_+^2\right).$$
We will refer to Newton diagrams simply as diagrams. Since a Newton diagram is uniquely determined by the compact faces of its border, we will often refer only to these compact faces.

We say that a diagram $\Gamma$ is supported by a set $\mathcal{S}$ if $\Gamma$ is the smallest diagram containing every point $P\in\mathcal{S}$. We say that $\Gamma$ lies below $\Sigma$ if $\Sigma\subset \Gamma$.

Let us denote by $(P_1,\dots,P_n)$ a diagram supported by points $P_1,\dots,P_n$. If $\Gamma$ is a diagram we will write $\Gamma+(P_1,\dots,P_n)$ for a diagram supported by $\supp \Gamma \cup \{P_1,\dots,P_n\}$. Any such diagram will be called a deformation of the diagram $\Gamma$.

If $P=(p,0)$, $Q=(0,q)$ then any translation of the segment $P,Q$ will be denoted as  $\tr(p,q)$, in other words
\begin{eqnarray}
\tr(p,q)& := &\text{hypotenuse of a right triangle with base}\nonumber\\
 &  &\text{of length }p\text{ and heigth }q\nonumber
\end{eqnarray}
We will write $n\tr(p,q)$ instead of $\tr(np,nq)$. Moreover, for $\tr(p_1,q_1),\dots,\tr(p_l,q_l)$ denote by 
$$(-1)^k\left(\ \tr(p_1,q_1)+\dots+\tr(p_l,q_l)\ \right)$$
any translation of a polygonal chain with endpoints $Q$, $Q+(-1)^k[p_1,-q_1]$, $\dots$ , $Q+(-1)^k\left[\sum_{i=1}^l p_i\ ,\ -\sum_{i=1}^l q_i\right].$

Note that if $(-1)^k=1$ we list the segments from top to bottom and if the sequence of the slopes $q_i/p_i$ is increasing, then $\tr(p_1,q_1)+\dots+\tr(p_l,q_l)$ is a Newton diagram. We will also write $\tr(P,Q)$ instead of $\tr(p,q)$ when we want to indicate fixed endpoints $P$ and $Q$ of the segment $\tr(p,q)$.

\subsection{Newton diagrams of singularities}\label{subsecNewtdiag}

We say that $\Gamma$ is the Newton diagram of an isolated singularity $f(x,y)=\sum_{i,j} c_{ij}x^iy^j$ if $\Gamma$ is the diagram supported by the set $\supp f=\{P\in\Za^2:c_P\neq 0 \}$. In such a case denote it by $\Gamma(f)$. We will say that $f$ is nondegenerate if it is nondegenerate in the sense of Kouchnirenko, see \cite{Kush}. Note that \cite{WallNonDeg} defines nondegeneracy differently but the definitions are equivalent in dimension $2$, see \cite{GreuelNonDegN2}.

A Newton diagram of a singularity is at distance at most $1$ from any axis.

\subsection{Newton numbers}\label{subsecNewtonNumber}

For a diagram $\Gamma\subset\re^2_+$, such that it has common points with both axis, its Newton number $\nu(\Gamma)$ is equal to
$$2A-p-q+1,$$
where $A$ is the area of the compliment $\re^2_+\setminus \Gamma$ and $p,q$ are the non-zero coordinates of the points of intersection. 

For any diagram $\Gamma\subset\re^2_+$ let $\nu_{p,q}(\Gamma)$ be the Newton number of a Newton diagram of $\Gamma+\{ (p,0), (0,q) \}$. Note that if $\Gamma$ is a diagram of an isolated singularity, then the definition does not depend on the choice of $p$ or $q$ if they are large enough, see \cite{Len08}. Hence the Newton number $\nu(\Gamma)=\nu_{p,q}(\Gamma)$, where $p,q$ sufficiently large, is well defined in the general case.

The motivation to study Newton numbers was given by Kouchnirenko in \cite{Kush}. In particular,
\begin{fact}\label{factNewMil}
For an isolated nondegenerate singularity the Newton number of its diagram and its Milnor number are equal.
\end{fact}

Similarly as for Milnor numbers, for a diagram $\Gamma$ consider the strictly decreasing sequence $(\nu_i)_{i=0,\dots,s}$ of positive integers attained as Newton numbers of deformations of $\Gamma$. Of course, $\nu_0=\nu(\Gamma)$ and $\nu_s=1$. The sequence $(\nu_{i-1}-\nu_{i})_{i=1,\dots,s}$ is the sequence of minimal jumps of Newton numbers.

Now for two useful properties.
\begin{prty}
\begin{enumerate}\label{prtyNewtonNr}\label{prty33}
\item 
If $\Sigma$ lies beneath $\Gamma$, then for any system of points $P_1,\dots,P_n$ the diagram $\Sigma+(P_1,\dots,P_n)$ lies below $\Gamma+(P_1,\dots,P_n)$ and the diagrams have common endpoints provided $\Sigma$ and $\Gamma$ had common endpoints.
\item 
If $\Sigma$ lies below $\Gamma$ and they have common endpoints, then the difference of Newton numbers $\nu_{p,q}(\Gamma) - \nu_{p,q}(\Sigma)$ is twice the difference of their areas. 
\end{enumerate}
\end{prty}


\subsection{General combinatorial remarks}\label{subsecGeneralComb}

Since a Newton number can be computed from the diagram, we will give some classic combinatorial tools that will help us in doing so.

\begin{fact}[Pick's Formula]
The area of a polygon with vertices from the lattice $\Za^2$ is equal to
$${B\over 2}+W-1,$$
where $B$ is equal to the number of points of the lattice $\Za^2$ which lie on its border and $W$ is the number of points of the lattice $\Za^2$ which lie in the interior.
\end{fact}

\begin{rk}[Tile Argument]
Consider a rhomboid $R(p,q)$ with vertices $(p,0)$, $(p-a,b)$, $(0,q)$, $(a,q-b)$, where $bp-aq=\pm 1$. The family 
$$\mathcal{R}(p,q)= \{\ R(p,q)+i[a,-b]+j[-(p-a),q-b] : \ i,j\in\Za\ \}$$
covers the real plane 
and consists of rhomboids with pairwise disjoint interiors. Moreover, every point in $\Za^2$ is a vertex of some rhomboid from this family.
\end{rk}
Indeed, since the area of $R(p,q)$ is $|pq-bp-(p-a)q|=1$, Pick's Formula implies that $R(p,q)\cap \Za^2$ is equal to the set of four vertices of $R(p,q)$. The rest follows immediately.

\subsection{EEA}\label{subsecEEA}
Let us recall the Extended Euclid's Algorithm. Note that $p,q$ being coprime implies $p,q\neq 1$.

\begin{fact}[Extended Euclid's Algorithm]
\label{factClassicEEA}
Take positive integers $p$ and $q$ which are coprime and $q>p$. The EEA goes as follows
$$
\begin{array}{r|ccccccc}
\text{variables} & P & Q & A' & A & B' & B & N \\\hline
\text{initial condition} & p & q & 0 & 1 & 1& 0 & \left\lfloor {q\over p} \right\rfloor \\
\\
\text{as long as $P\neq 0$ substitute} & Q-NP & P & A-NA' & A' & B-NB' & B' &  \left\lfloor {P\over Q-NP} \right\rfloor \\
\\
\text{the output line } P=0 & 0 & 1 & \pm\ p & \mp\  a & \mp\ q & \pm\  b & 0 
\end{array}
$$
Positive integers $a,b$ in the last line are such that $a<p, b<q$ and $|bp-aq|=1$.
\end{fact}

We will adjust the algorithm to our needs. Reverse the order of the lines and number them from $0$ for the output line to $k_0+2$ for the initial conditions line (we always have at least 3 lines, hence $k_0\ge 0$). Note that $a_0=p, b_0=q$ and we get a modified table
\begin{equation}
\label{EEAFullGeneral}
\begin{array}{cc|l}
p & q & \\\hline
a_1 & b_1 & n_1\\
\vdots & \vdots & \vdots\\
a_{k_0+1} & b_{k_0+1} & n_{k_0+1}\\
a_{k_0+2} & b_{k_0+2} & 
\end{array}
\end{equation}
which consists of columns $A', B'$ and $N$ from original EEA in reverse order and dropping the signs. Note that $a_1=a$ and $b_1=b$. 

Consider an example that we will use as an illustration throughout this paper.
\begin{ex}\label{przyEEA} For $p=40$ and $q=73$ we have $k_0=4$ and
$$
\begin{array}{cc|l}
40 & 73 & \\\hline
17 & 31 & 2\\
6 & 11 & 2\\
5 & 9 & 1\\
1 & 2 & 4\\
1& 1 & 1\\
0&1
\end{array}
$$
In particular, $31\cdot 40-17\cdot 73=-1=(-1)^{4-0+1}$.
\end{ex}

We will list some properties of EEA adjusted to our notations.
\begin{prty}\label{wlEEA}
\begin{enumerate}
\item The values in the last two lines are always
$$
\begin{array}{lll}
a_{k_0+1}=1, & b_{k_0+1}=\left\lfloor {q\over p} \right\rfloor, & n_{k_0+1}=a_{k_0}\\
a_{k_0+2}=0, & b_{k_0+2}=1 &
\end{array}
$$
and necessarily $a_{k_0+1}b_{k_0+2}-b_{k_0+1}a_{k_0+2}=1$.
\item Each new line can be obtained as the rest from division from the former two lines (except $a_{k_0+1}$ and $b_{k_0+2}$).
In particular for any $k=1,\dots,k_0+1$ we have
$$
a_{k+1}=a_{k-1}-n_{k}a_{k},\quad b_{k+1}=b_{k-1}-n_{k}b_{k} ,
$$
$n_k= \left\lfloor {a_{k-1}\over a_k}\right\rfloor$ for $k< k_0$ and $n_k=\left\lfloor {b_{k-1}\over b_k} \right\rfloor$ for $k\le k_0$.
\item The positive integers $a_k$ and $b_k$ are coprime and the sign of $a_kb_{k+1}-b_ka_{k+1}$ alternates.
In particular, for $k=0,\dots,k_0$ we get $ b_{k}>a_k\ge 1$
and
$$a_kb_{k+1}-b_ka_{k+1}=(-1)^{k_0-k+1}$$
\end{enumerate}
\end{prty}

\subsection{Remarks on EEA}\label{subsecEEArks}

Let $q=mp+r$, where $r<p<q$ and $q,p$ are coprime. Consider EEA beginning with 
\begin{equation}
\label{EEA}
\begin{array}{cc|l}
p & q & \\\hline
a & b & n\\
a' & b' & n'\\
a'' & b'' & n''\\
\dots & \dots & \dots
\end{array}
\end{equation}
and denote $${\rm sign}(p,q):=bp-aq.$$
Note that ${\rm sign}(p,q)$ is equal $(-1)^{k_0+1}$. We will retain this notation throughout the rest of the paper and prove some technical properties that will be useful.
\begin{prty}
\label{prtyN}
We may assume that $n$ in EEA is >1.
\end{prty}
\noindent{\bf Proof.} Instead of \eqref{EEA} consider a shorter EEA
$$
\begin{array}{cc|l}
p & q & \\\hline
a' & b' & n'+1\\
a'' & b'' & n''\\
\dots & \dots & \dots
\end{array}
$$
We have $pb'-qa'= -n \cdot {\rm sign}(p,q) =-{\rm sign}(p,q)$ and $p=a+a'=n'a'+a'+a''$. Hence in the table above signs alternate and the table above has all properties listed in Property \ref{wlEEA}. The rest of the table does not change.\hfill$\blacksquare$

\begin{prty}
\label{prtyShortEEA}
EEA ends with
$$
\begin{array}{cc|l}
\dots & \dots & \dots\\
\tilde{a} & \tilde{a}m+1 & \tilde{n}\\
1 & m & \tilde{a}\\
0 & 1 & 
\end{array}
$$
where $\tilde{a},\tilde{n}$ are positive integers. Moreover, $\tilde{a}=1$ iff $2r>p$.
\end{prty}
\noindent{\bf Proof.} The last lines were already given in Property \ref{wlEEA}. We need to prove the second part. To study the last lines recall the classic EEA, Fact \ref{factClassicEEA}. It easily follows that $\tilde{a}=\lfloor p/r \rfloor$. Hence $\tilde{a}$ is equal to one if and only if $2r>p$.  \hfill$\blacksquare$

\begin{prty}
\label{prtyResztaRownaN}
\label{prtyShortEEAa1}
If $k_0=2$, EEA is of the form
\begin{equation}
\label{eqSchortEEA}
\begin{array}{cc|l}
p & q & \\\hline
a & am+1 & n=r\\
1 & m & n'=a\\
0 & 1 & 
\end{array}
\end{equation}
and $a=1$ iff $q=mp+p-1$.
\end{prty}
\noindent{\bf Proof.} Taking into account Property~\ref{prtyShortEEA} above one needs to show only that $n=r$ as well as the equivalence. Indeed, if $a''=0$, then by the above Property~\ref{prtyShortEEA} we get $\tilde{a}=a, \tilde{a}m+1=b$, $\tilde{n}=n$ thus $(am+1)p-aq=1$. Hence $1=a(mp-q)+p=-ra+p$ and it follows that $p=ra +1$ on one hand, while $p=na+1$ on the other. Moreover, if $a=1$, then $p(m+1)-q=1$. On the other hand, if $q=mp+p-1$ the EEA is of the form \eqref{eqSchortEEA} with $a=1$. \hfill$\blacksquare$

Moreover, as a special case of Property \ref{prtyShortEEA} we get
\begin{prty}
\label{prtyVeryShortEEA}
$q=mp+1$ if and only if EEA is of the form
\begin{equation}
\label{eqVerySchortEEA}
\begin{array}{cc|l}
p & q & \\\hline
1 & m & p\\
0 & 1 & 
\end{array}
\end{equation}
\end{prty}

Now for two technical properties
\begin{prty}
\label{prtyEEApjqj}
Take $j< n'$, any positive integer $l$, $p^j=l(a-ja')+a'$ and $q^j=l(b-jb')+b'$. Then $p^j,q^j$ are coprime and their EEA is of the form
$$
\begin{array}{cc|l}
l(a-ja')+a' & l(b-jb')+b' & \\\hline
a-ja' & b-jb' & l\\
a' & b' & n'-j\\
a'' & b'' & n''\\
\dots & \dots & \dots
\end{array}
$$
\end{prty}
Indeed, $(n(a-ja')+a')(b-jb') - (n(b-jb')+b')(a-ja') = a'(b-jb')-b'(a-ja')=a'b-b'a = {\rm sign}(p,q)$.

\begin{prty}
\label{prtyEEAatilde}
Take a positive integer $N$, assume $a''\neq 0$. Then $Na'+a''$ and $Nb'+b''$ are coprime and their EEA is of the form
$$
\begin{array}{cc|l}
Na'+a'' & Nb'+b'' & \\\hline
a' & b' & N\\
a'' & b'' & n''\\
\dots & \dots & \dots
\end{array}
$$
\end{prty}
Indeed, $(Na'+a'')b''-(Nb'-b'')a''=a''b'-a'b''=-{\rm sign}(a',b')$.

Now it is easy to see that if $a\neq 1$, then $a'\neq 0$.

\section{Main steps of proof}\label{secLemmas}

Choose the line $i\le k_0$ in the EEA \eqref{EEAFullGeneral} for $a_0,b_0$, where $1<a_0<b_0$ are coprime. Denote $p=a_i,q=b_i$ and 
assume EEA is of the form \eqref{EEA}. We have $bp-aq=(-1)^{k_0-i+1}$, recall ${\rm sign}(p,q)=bp-aq$. We will consider deformations of $\tr(p,q)$.

Let $Q$ denote the upper and $P$ denote the lower endpoint of the diagram $\tr(p,q)$ if ${\rm sign}(p,q)=-1$, reversely if ${\rm sign}(p,q)=1$.

\subsection{Decreasing $p$ and $q$} \label{subsecDecrease}  
In this paragraph, informally speaking, we will aim at replacing $p=na+a'$ by $p=n(a-a')+a'$ (and at the same time $q=nb+b'$ by $n(b-b')+b'$). In the next paragraph~\ref{subsecReductionLine} we will prove that one can do it recursively until $a-ka'=a''$. This will allow us to use EEA and reduce the problem to repetition of the procedure for consecutive levels of the EEA table~\eqref{EEAFullGeneral}.

Consider a diagram
$$ \Gamma^k = -{\rm sign}(p,q)\big( k\tr(a',b')+\tr(p-ka,q-kb)+k\tr(a-a',b-b')\big)$$
where $0\leq k\leq n$. Denote also the points
$$P^k=P-{\rm sign}(p,q) k[-(a-a'),b-b'],\quad Q^k=Q-{\rm sign}(p,q)k[a',-b'], $$
in the support of $\Gamma^k$ such that
$$\Gamma^k = -{\rm sign}(p,q)\left(\tr(Q,Q^k)+\tr(Q^k,P^k)+\tr(P^k,P)\right).$$
Note that $\Gamma^0=\tr(p,q)$, 
$$\Gamma^n=-{\rm sign}(p,q)\big(\ (n+1)\tr(a',b')+n\tr(a-a',b-b')\ \big)$$
and every $\Gamma^k$ is a Newton diagram.

Consider points
$$P_i^k = P^k-{\rm sign}(p,q)\cdot i[-a,b],\quad i=1,\dots,n-k$$
and
$$D_i^k = P_i^k+{\rm sign}(p,q)[-a',b'],\quad i=1,\dots,n-k.$$

\begin{proc}\label{procka1}
Consecutively for $k=0,\dots,n-1$ take diagrams
$$\Gamma^k+P_i^{k}\quad i=1,\dots,n-k,$$
$$\Gamma^k+D_i^{k}\quad i=1,\dots,n-k.$$
\end{proc}

Note that $$\Gamma^k+P_i^{k} = \Gamma+(P_{n-k-1}^{k-1},P_i^{k},D_1^{k-1}) = \Gamma + (P^k,P^k_i,Q^k)$$ and analogously $$\Gamma^k+D_i^{k} = \Gamma+(P_{n-k-1}^{k-1},D_i^{k},D_1^{k-1}) = \Gamma + (P^k,D^k_i,Q^k).$$

\begin{prop}\label{propReductioAtoA1}
If $a\neq 1$ the choice of deformations in Procedure \ref{procka1} gives the opening terms of the sequence of minimal jumps of Newton numbers
$$\underbrace{1\ ,\ \dots\ ,\ 1}_{n(n+1)}$$
\end{prop}

Proof will follow after some lemmas below.

\begin{lem}\label{lemPropertiesP}
 For any fixed $k$ we have
\begin{enumerate}
\item For $i=1,\dots,n-k-1$ the deformation $\Gamma^k+P_i^k$ has the diagram
$$-{\rm sign}(p,q)(\tr(Q,Q^k)+\tr(Q^k,P_i^k)+\tr(P_i^k,P^k)+\tr(P^k,P)).$$
\item The deformation $\Gamma^k+P_{n-k}^k$ has the diagram
$$-{\rm sign}(p,q)\big((k+1)\tr(a',b')+(n-k)\tr(a,b) +k\tr(a-a',b-b')\big).$$
\end{enumerate}
\end{lem}

\noindent{\bf Proof.} First note that $P_1^k,\dots,P_{n-k}^k$ and $P^k$ are colinear. Moreover, from Euclid's Algorithm $p-ka=(n-k)a + a'$ and $q-kb=(n-k)b+b'$ with $a'b'\neq 0$ and $k=0,\dots,n$. To prove (1) it suffices to note that as a consequence the slopes of $\tr(P^k,P), \tr(P^k,P_i^k), \tr(P_i^k,Q^k)$ and $\tr(Q^k,Q)$ exactly in that order constitute a strictly monotone sequence. Point (2) follows from the above considerations taking into account the fact that $p-na=a'$ and $q-nb=b'$, hence the slopes of $\tr(P_{n-k}^k,Q^k)$ and $\tr(Q^k,Q)$ are equal.\hfill$\blacksquare$\\

\begin{lem} \label{lemjumpP}
For fixed $k$ and $i=1,\dots,n-k$ we have 
$$\nu(\Gamma^k)-\nu(\Gamma^k+P_i^k)=i.$$
\end{lem}

\noindent{\bf Proof.} Note that from Lemma \ref{lemPropertiesP} it follows that we add only points that are in the interior of the triangle with hypotenuse $\tr(p-ka,q-kb)$. Moreover, they all lie on or over the line passing through $Q^k$ with the slope as of $\tr(a',b')$ and on or over the line passing through $P^k$ with the slope as in $\tr(a-a',b-b')$. Hence the difference of Newton numbers of $\Gamma^k$ and $\Gamma^k+P_i^k$ is equal to double the area of their difference. 

Now the claim easily follows from Tile Argument and Pick's formula, since double the area of the triangle $P^kQ^kP^k_{n-k}$ is equal $n-k$.\hfill$\blacksquare$\\

\begin{lem}\label{lemPropertiesD}
For any fixed $k$ we have
\begin{enumerate}
\item For $i=2,\dots,n-k$ the deformation $\Gamma^k+D_i^k$ has the diagram
$$-{\rm sign}(p,q)(\tr(Q,Q^k)+\tr(Q^k,D_i^k)+\tr(D_i^k,P^k)+\tr(P^k,P)).$$
\item The above holds for $i=1$ provided $a\neq 1$.
\end{enumerate}
\end{lem}

\noindent{\bf Proof.} First note that from their definition, the points $D_1^k, \dots, D_{n-k}^k$ all lie on a translation of the segment with endpoints $P_1^k, P_{n-k}^k$ by the vector $[-a',b']$. Hence to prove (1) it suffices to note that $Q,Q^k$ and $D_{n-k}^k$ are colinear and the slope of $\tr(D_2^k,P^k)$ is bigger then that of $\tr(P^k,P)$ in the case ${\rm sign}(p,q)=-1$ (smaller in the other case).

Note that we have $a\neq a'$ unless $a=a'=1$. Hence if $a\neq a'$, we have equality of the slopes of $\tr(D_{1}^k,P^k)$ and $\tr(P^k,P)$. Which proves (2).\hfill$\blacksquare$\\

\begin{rk} If $a\neq 1$
$$\Gamma^k+(P_{n-k}^k, D_1^k)=\Gamma^0 + (P_{n-k}^k, D_1^k) = \Gamma^{k+1}$$
and $\Gamma^n$ lies below all points considered above.
\end{rk}

\begin{lem} \label{lemjumpD}
For fixed $k$ if $a\neq 1$, then for $i=1,\dots,n-k$ we have 
$$\nu(\Gamma^k)-\nu(\Gamma^k+D_i^k)=n-k+i.$$
\end{lem}

\noindent{\bf Proof.} Similarly to the opening argument of the proof of Lemma \ref{lemjumpP} we derive from Lemma \ref{lemPropertiesD} that $\nu(\Gamma^k)-\nu(\Gamma^k+D_i^k)$ is equal to double the area of the difference of the diagrams. Moreover, this difference can be computed when considering only the segment $\tr(p-ka,q-kb)$.

Consider double the area of $P^kQ^kD^k_{j}$ with fixed $j$. We will compute it using Pick's formula. Without loss of generality we can assume that $bp-aq=-1$. 

Note that due to Tile Argument, the only points that may lie in the triangle $P^kQ^kD^k_{j}$ are the points $P_i^k$. First, note that any $P_i^k$ with $i<j$ lies in the interior of the triangle with vertices $P^kQ^kD^k_{j}$, since it suffices to notice that the segments $P^kP^k_j$ and $D_i^kD_j^k$ are parallel. Any $P_i^k$ with $i\ge j$ lies in the interior of the triangle with vertices $P^kQ^kD^k_{j}$ if and only if the slope of $Q^kD^k_j$ is greater than the slope of $Q^kP^k_i$ (in absolute values) i.e.
\begin{equation}
\label{EQlemma4slopes}
(q-kb-jb+b')(p-ka-ia)>(p-ka-ja+a')(q-kb-ib),
\end{equation}
whereas $P_i^k$ lies on its side if and only if there is an equality of the slopes. Equation~\eqref{EQlemma4slopes} is equivalent to 
$$(i-j)(b(p-ka)-(q-kb)a)+b'(p-ka)-a'(q-kb)-i>0.$$
Note that from $q=nb+b'$ and $p=na+a'$ it follows that $a'q-b'p=-n$. Hence $a'(q-kb)-b'(p-ka)=k-n$ and \eqref{EQlemma4slopes} is equivalent to 
$${j+n-k\over 2}>i. $$
Of course, equality in \eqref{EQlemma4slopes} holds if and only if ${j+n-k\over 2}=i$. By $\#$ denote the number of elements. Above combined with Pick's formula gives that double the area of the triangle with vertices $P^kQ^kD^k_{j}$ is
$$B+2W-2 = 3 + \#\left\{i\ |\  {j+n-k\over 2}=i\ge j \right\}+ 2(j-1)+2\#\left\{i\ |\  {j+n-k\over 2}>i\ge j \right\}  -2.$$
If $j+n-k$ is even, then the above is equal to 
$$1+ 1+ 2(j-1)+2\left({j+n-k\over 2}-1-j+1\right)=j+n-k.$$
If $j+n-k$ is odd, then the above is equal to 
$$1+ 0 + 2(j-1)+2\left({j+n-k-1\over 2}-j+1\right)=j+n-k.$$
This gives the assertion.\hfill$\blacksquare$\\

\begin{rk}\label{rkJumps}
In particular, for $a\neq 1$ Lemmas \ref{lemjumpP} and \ref{lemjumpD} imply that the sequence of minimal jumps for diagram $\Gamma^k$ begins with
$$\underbrace{\ 1\ ,\ \dots\dots\ ,\ 1\ }_{2(n-k)}.$$
\end{rk}

\noindent{\bf Proof of Proposition \ref{propReductioAtoA1}.} Thanks to Lemmas \ref{lemjumpP} and \ref{lemjumpD} (see Remark \ref{rkJumps}) we only have to show that 
$$\nu(\Gamma_k)-\nu(\Gamma_{k+1})=2(n-k)$$
for $k=0,\dots,n-1$. We compute this number as double the area of the polygon with vertices $P_{n-k}^{k}, D_1^{k}, P_{n-k-1}^{k+1}, D_1^{k+1}$. From Tile Argument for $[a,-b]$ it follows that the only integer points on the boundary are the vertices, whereas $P_k^j$ for $j=1,\dots,n-k-1$ lie in the interior. Again from Tile Argument for $[a',-b']$ these points are the only integer ones to lie there. Therefore, from Pick's formula we get 
$$\nu(\Gamma_k)-\nu(\Gamma_{k+1}) = 4+2(n-k-1) -2=2(n-k).$$
Therefore consecutive choices in Procedure \ref{procka1} give consecutively $1$ in the sequence of minimal jumps and
$$\nu(\Gamma^0) - \nu(\Gamma^n)=\sum_{k=0}^{n-1}2(n-k)=n(n+1).$$ 
This ends the proof.\hfill$\blacksquare$\\

\subsection{Reduction of the line in EEA} \label{subsecReductionLine}
We will now recursively substitute $p$ by $n(a-a')+a'$ (compare previous Subsection~\ref{subsecDecrease}) i.e. we will reduce $p$ to $a$. 

Consider diagrams
$$\Sigma^j=-{\rm sign}(p,q)\big(nj\ \tr(a',b')+\tr\left(\ a'+n(a-ja'),\ b'+n(b-jb')\ \right)\big)$$
for $j=0,\dots,n'$. Let 
$$P(\Sigma^j)=Q-{\rm sign}(p,q)nj[a',-b']$$ 
be the point such that
$$\Sigma^j=-{\rm sign}(p,q)(\tr\left(Q,P(\Sigma^j)\right)+\tr\left(P(\Sigma^j),P\right)).$$

Note that $\Sigma^0=\tr(p,q)$, $$\Sigma^1=-{\rm sign}(p,q)(n\ \tr(a',b')+\tr\left(\ a'+n(a-a'),\ b'+n(b-b')\ \right))$$ and $$\Sigma^{n'}=-{\rm sign}(p,q)(nn'\tr(a',b')+\tr(a'+na'',b'+nb'')).$$
Every $\Sigma^j$ is a diagram.

\begin{figure}[h]
\centering
\includegraphics [height=6cm]{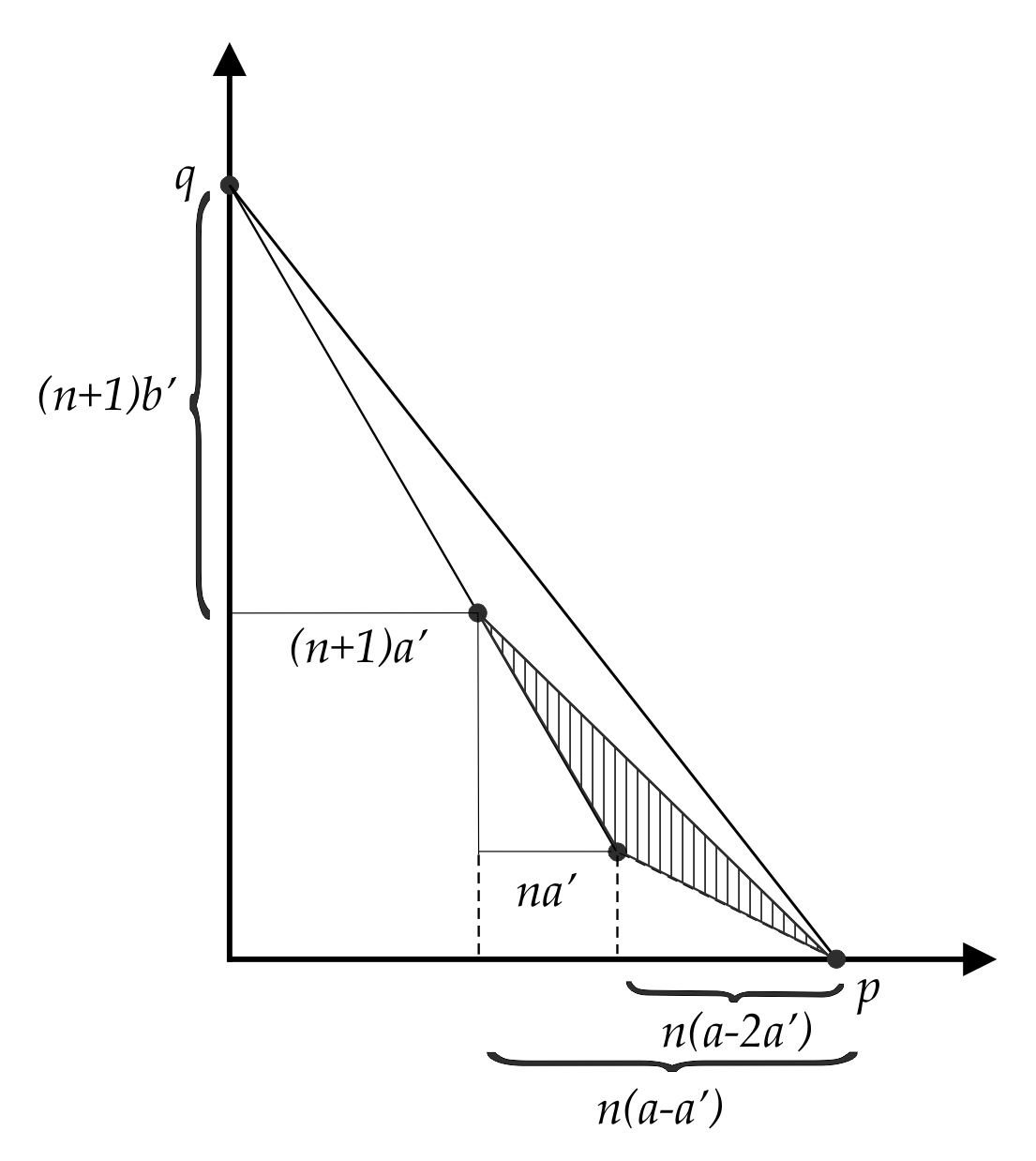} \label{r:1}
\caption{Procedure \ref{procka2}, a step from $\Sigma^1$ to $\Gamma^n(\Sigma^1)$. It is simply an iteration of Procedure~1. The shaded area is the difference between the two diagrams.}\label{f:1}
\end{figure}

Let $p^j=a'+n(a-ja')$, $q^j=b'+n(b-jb')$ and
$$
\begin{array}{c}
\Gamma^k(\Sigma^j) =-{\rm sign}(p,q)   \Big( (nj + k)\tr(a',b')  + \tr\left(p^j-k(a-ja'),q^j-k(b-jb')\right) +\\ 
  + k\tr\left(a-(j+1)a',b-(j+1)b'\right) \Big).
\end{array}
$$
Hence
$$\Gamma^k(\Sigma^j) = -{\rm sign}(p,q)( \tr(Q,Q^k(\Sigma^j))+\tr(Q^k(\Sigma^j),P^k(\Sigma^j))+\tr(P^k(\Sigma^j),P) ),$$
where
$$P^k(\Sigma^j)=P-{\rm sign}(p,q)\cdot k[-(a-(j+1)a'),b-(j+1)b'],$$  
$$Q^k(\Sigma^j)=Q-{\rm sign}(p,q)\cdot (nj+k)[a',-b'].$$

Consider points
$$P_i^k(\Sigma^j) = P^k(\Sigma^j)-{\rm sign}(p,q)\cdot i[-(a-ja'),b-jb'],\quad i=1,\dots,n-k$$
and
$$D_i^k (\Sigma^j)= P_i^k(\Sigma^j)+{\rm sign}(p,q)[-a',b'],\quad i=1,\dots,n-k.$$

\begin{proc}\label{procka2}
Fix $j\in\{0,\dots,n'-1\}$. Consecutively for $k=0,\dots,n-1$ take diagrams
$$\Gamma^k(\Sigma^j)+P_i^{k}(\Sigma^j)\quad i=1,\dots,n-k,$$
$$\Gamma^k(\Sigma^j)+D_i^{k}(\Sigma^j)\quad i=1,\dots,n-k.$$
\end{proc}

Note that 
$$\Gamma^n(\Sigma^j)= -{\rm sign}(p,q)\big( \left(n(j+1)+1 \right) \tr (a',b') + n\ \tr\left(a-(j+1)a',b-(j+1)b'\right) \big).$$

\begin{rk}
\label{rkAllDeformationsAboveGammaN}
All points $P_i^{k}(\Sigma^j)$ and $D_i^{k}(\Sigma^j)$ lie on or above the diagram $\Gamma^n(\Sigma^{n'-1})= -{\rm sign}(p,q) ( \left(nn'+1 \right) \tr (a',b') + n\ \tr\left(a'',b''\right) ).$
\end{rk}

Below is a generalisation of Proposition \ref{propReductioAtoA1}.

\begin{prop}
\label{propSigmaKdoSigmaK1}
If $a\neq 1$ the choice of deformations in Procedure \ref{procka2} for the diagram $\Sigma^j$ gives the opening terms of the sequence of minimal jumps of Newton numbers
$$\underbrace{1\ ,\ \dots\ ,\ 1}_{n(n+1)}$$
provided $j<n'-2$ or $j=n'-1$ and $a''\neq 0$.
\end{prop}
\noindent{\bf Proof.} Apply Proposition \ref{propReductioAtoA1} to $\tr\left( p^j,q^j \right)$, where $$p^j=a'+n(a-ja')\quad {\rm and}\quad q^j=b'+n(b-jb').$$ From Property \ref{prtyEEApjqj} the last diagram $\Gamma^n$ is of the form
$$-{\rm sign}(p,q) \big((n+1)\tr(a',b')+n\tr(a-ja'-a',b-jb'-b') \big).$$
Hence $nj\tr(a',b')+\Gamma^n$ is a diagram. Moreover, it is exactly $\Gamma^n(\Sigma^j)$ and no segment lies on any axis if $j<n'-2$ or $j=n'-1$ and $a''\neq 0$. Therefore, all preceding $\Gamma^k$ for $k=0,\dots, n-1$ coupled with $nj\tr(a',b')$ are also diagrams (in fact equal to $\Gamma^k(\Sigma^j)$). Hence the claim follows from Property \ref{prtyNewtonNr} and Proposition \ref{propReductioAtoA1}.\hfill$\blacksquare$\\

\begin{prop}
\label{propSigma0DoSigmaN}
Let $a\neq 1$ and $a''\neq 0$. For consecutive $j=0,\dots,n'-1$ consider points as in Procedure \ref{procka2}. They give the opening terms for $\tr(p,q)$ of the sequence of minimal jumps of Newton numbers
$$\underbrace{1\ ,\ \dots\ ,\ 1}_{n(nn'+1)}$$
\end{prop}
Proof will follow immediately from 

\begin{lem}
\label{lemSigmaGamma}
For $j=0,\dots,n'-1$ the diagram $\Gamma^n(\Sigma^j)$ lies below $\Sigma^{j+1}$. Moreover, if $j<n'-2$ or $j=n'-1$ and $a''\neq 0$ we have 
$$\nu(\Sigma^{j+1})=\nu(\Gamma^n(\Sigma^j))+n.$$ 
\end{lem}
\noindent{\bf Proof.} From Property \ref{prtyNewtonNr} and the form of the diagrams one has to compute double the area of the triangle with vertices $P(\Sigma^{j+1}), P(\Sigma^{j+1})-{\rm sign}(p,q)[a',-b']$ and $Q$. From Tile Argument for $p^{j+1}$ and $q^{j+1}$ as well as for $a'$ and $b'$, the only points that lie in this triangle lie on its sides and there are exactly $n+2$ such points. Hence form Pick's formula we get the claim.\hfill$\blacksquare$\\

\noindent{\bf Proof of Proposition \ref{propSigma0DoSigmaN}.} The claim follows from Lemma \ref{lemSigmaGamma} above and the fact that we have $n'$ steps. Each step gives $n^2+n$ ones in the sequence (see Proposition \ref{propSigma0DoSigmaN}), where $n$ ones are attained twice with the last deformations of $\Sigma^j$ and initial deformations of $\Sigma^{j+1}$ (see Lemma \ref{lemSigmaGamma} above) with the exception of the $(n'-1)$th step. Since $a''\neq 0$, none of the points lie on an axis.\hfill$\blacksquare$\\

\begin{figure}[h]
\centering
\includegraphics [height=6cm]{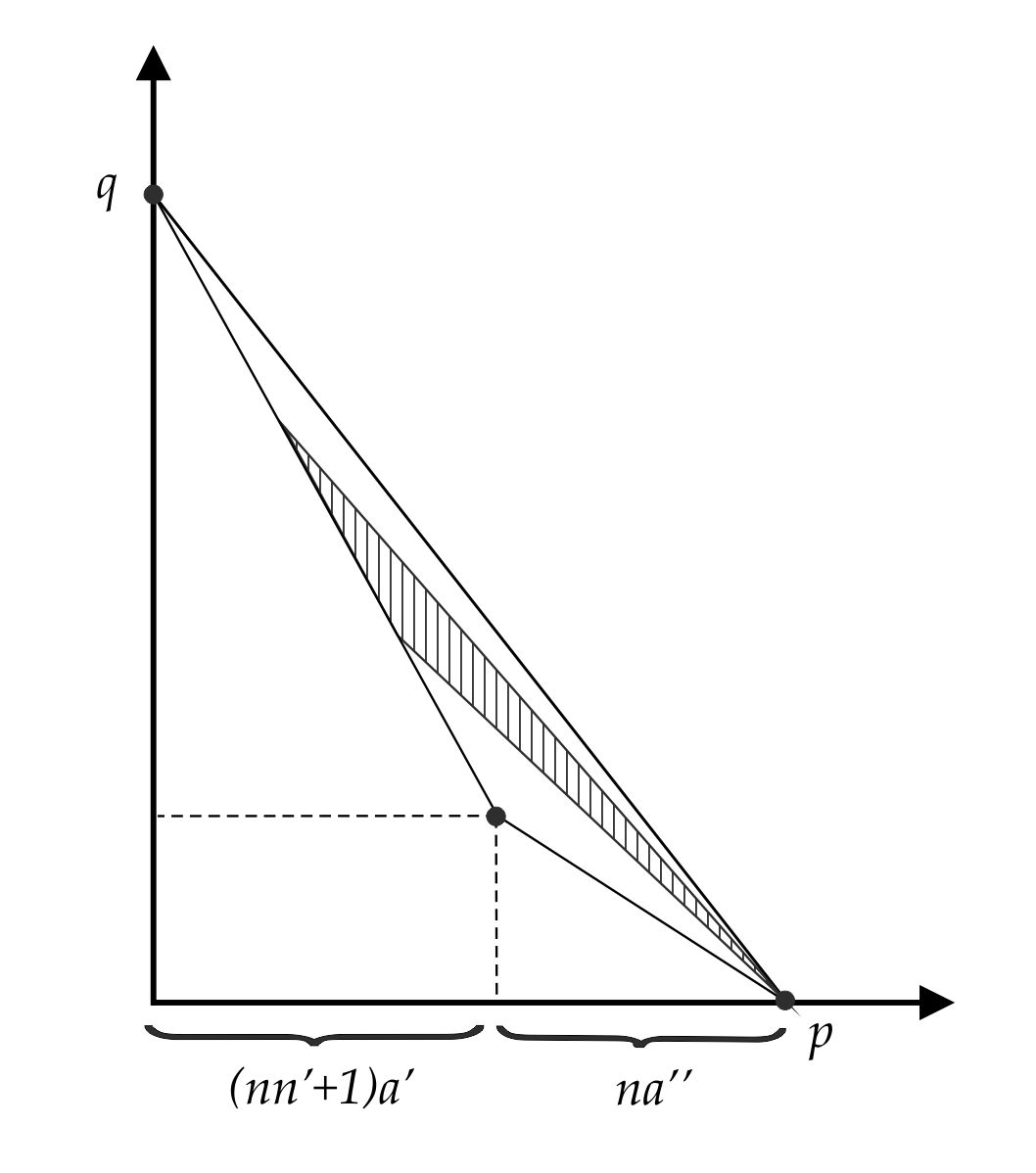} \label{r:2}
\caption{Procedure \ref{procka2} ends with the diagram $\Sigma^{n'}$ provided $a''\neq 0$. The shaded area represents a difference between the diagrams in intermediate steps, compare Figure \ref{r:1}.}
\end{figure}

Suppose $a''\neq 0$. Consider the diagram
$$\Theta=-{\rm sign}(p,q) \big(\tr\left((nn'+1)a'+a'', (nn'+1)b'+b'' \right) + (n-1) \tr(a'',b'') \big). $$

Note that $n-1>0$ from Property \ref{prtyN}. Moreover, from Property \ref{prtyEEAatilde} we get that EEA is the same as EEA of $a,b$ up to the first line.

\begin{lem} 
\label{lemTheta}
If 
$a''\neq 0$, then $\Theta$ lies above $\Gamma^n(\Sigma^{n'-1})$ and
$$\nu(\Theta)-\nu(\Gamma^n(\Sigma^{n'-1}))=nn'+1$$
\end{lem}
\noindent{\bf Proof.} Note that $(nn'+1)a'+a''$ and $(nn'+1)b'+b''$ are coprime and their EEA series is given in Property \ref{prtyEEAatilde}. Hence using Tile Argument and Pick's formula we easily get the claim.\hfill$\blacksquare$\\


\subsection{Short EEA} \label{subsecShortEEA} 
In this section we consider two cases left i.e. what happens if $a''=0$ or $a'=0$. Note that $a=1$ implies $a''=0$ or $a'=0$ i.e. EEA is of the form \eqref{eqSchortEEA} or \eqref{eqVerySchortEEA}.

Let us remind that $q=mp+r$, where $0<r<p$. From Property \ref{prtyResztaRownaN} we have $r=n$ for short EEA. Under notation of Procedure \ref{procka1} consider
\begin{proc}
\label{prockaShortNie1}
Let $a''=0$ and $ a\neq 1$. For $j\in\{0,\dots,n'-2\}$ consecutively for $k=0,\dots,n-1$ take diagrams
$$\Gamma^k(\Sigma^j)+P_i^{k}(\Sigma^j)\quad i=1,\dots,n-k,$$
$$\Gamma^k(\Sigma^j)+D_i^{k}(\Sigma^j)\quad i=1,\dots,n-k.$$
\end{proc}

Note that the last diagram in the procedure above is $\Gamma^n(\Sigma^{n'-2}) $ of the form 
\begin{equation}
\label{eqDiagramKoncowy}
n\tr(1,m+1)+(n(n'-1)+1)\tr(1,m).
\end{equation}

\begin{proc}
\label{prockaShortA1}
Let $a'\neq 0$ and $ a= 1$. Consecutively take diagrams
$$(P,Q)+Q_i\quad i=1,\dots,p-1,$$
where $Q_i=Q-i[-1,m+1].$
\end{proc}

Note that $a'=0$ iff $q=mp+1$.

\begin{proc}
\label{prockaVeryShortEEA}
Let $a'=0$. Take diagrams
$$(P,Q)+P_i \quad i=1,\dots,p-1,$$
where $P_i=P+i[-1,m].$
\end{proc}

\begin{figure}[h]
\centering
\includegraphics [height=6cm]{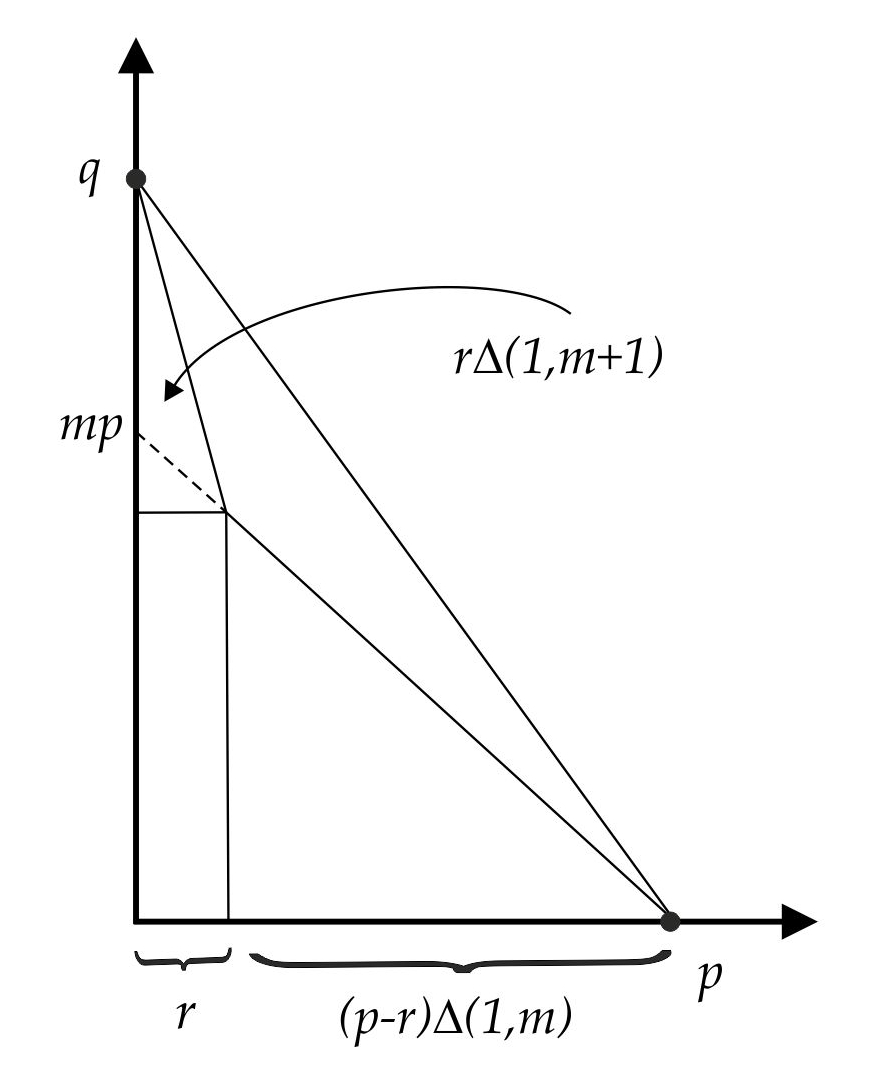} \label{r:4}
\caption{Procedures~\ref{prockaShortA1}, \ref{prockaShortNie1} and \ref{prockaVeryShortEEA} all end with this diagram, where $q=mp+r$.}
\end{figure}

\begin{prop}
\label{propShortEEA}
If $a''=0$ or $a'=0$, the opening terms of the sequence of minimal jumps of Newton numbers are
$$\underbrace{1\ ,\ \dots\ ,\ 1}_{r(p-r)}$$
\end{prop}

\noindent{\bf Proof.} We have three cases.

If $q=kp+1$ i.e. $a'=0$ consider the choice of deformations in Procedure \ref{prockaVeryShortEEA} and the claim follows immediately from Lemma \ref{lemjumpP} and Property \ref{prtyVeryShortEEA}. Note that the number of jumps above is also equal to $r(p-r)$, because here $r=1$ and the diagram $(P,Q)+P_{p-1}$ is of the form \eqref{eqDiagramKoncowy}.

If $a'\neq 0$ and $ a= 1$, then $a'=1, a''=0$. Consider the choice of deformations in Procedure \ref{prockaShortA1} and the claim follows immediately from Lemma \ref{lemjumpP}. Note that the number of jumps above is also equal to $r(p-r)$, because here $r=p-1$, see Property~\ref{prtyShortEEAa1}. Again, the diagram $(P,Q)+Q_{p-1}$ is of the form \eqref{eqDiagramKoncowy}.

If $a''=0$ and $ a\neq 1$ consider deformations in Procedure \ref{prockaShortNie1}, the proof is the same as in Proposition \ref{propSigma0DoSigmaN} and follows from Lemma \ref{lemSigmaGamma}. The length of the sequence of jumps is hence equal to $(n'-1)n^2+n=n(na-n+1)=r(p-r)$ thanks to Properties~\ref{prtyResztaRownaN} and \ref{prtyShortEEA}.

Hence the claim. \hfill$\blacksquare$\\

\begin{rk}
\label{rkLastDiagramShort}
Note that if $a''=0$ or $a'=0$, the last diagram is of the form \eqref{eqDiagramKoncowy}.
\end{rk}

\section{Main theorem combinatorially}\label{secMainThmComb}

\begin{thm}
\label{thmMainCombinatorially}
Given a diagram $\tr(a_0,b_0)$, where $a_0,b_0$ are coprime, the sequence of minimal jumps of Newton numbers commences with 
$$\underbrace{1\ ,\ \dots\ ,\ 1}_{r(a_0-r)}$$
where $r$ is the rest out of division of $b_0$ by $a_0$.
\end{thm}
\noindent{\bf Proof.} Suppose EEA is of the form \eqref{EEAFullGeneral}. Consider an auxiliary sequence
$$z_0=1,\quad, z_1=n_1,\quad, z_k=z_{k-2} + z_{k-1}n_k .$$
This sequence coincides with the column $P$ in reverse order in EEA, see Fact \ref{factClassicEEA}. Note that by Property \ref{prtyN} we may assume that $z_1>1$ and hence $(z_k)$ is strictly increasing.
\begin{proc}
\label{prockaWszystko}
Take $k=1$. 

Put $L=z_{k-1}-1$,
$$
\begin{array}{ll|l}
p=z_ka_k+a_{k+1} & q=z_kb_k+b_{k+1} & \\\hline
a=a_k & b=b_k & n=z_k\\
a'=a_{k+1} & b'=b_{k+1} & n'=n_{k+1}\\
a''=a_{k+2} & b''=b_{k+2} &
\end{array}
$$
and consider deformations of the diagrams 
$$\Theta_k=-{\rm sign}(p,q) \big( L\tr(a',b') + \tr(p,q) \big).$$

For $\tr(p,q)$:

If $a'=0$, use Procedure \ref{prockaVeryShortEEA}.

If $a''=0$, $a'=1$ and $a=1$ use Procedure \ref{prockaShortA1}.

If $a''=0$, $a'=1$ and $a\neq 1$ use Procedure \ref{prockaShortNie1}.

Otherwise, use Procedure \ref{procka2}, afterwards substitute $k$ by $k+1$ and proceed as above.
\end{proc}
This procedure will end, because the EEA sequence is finite. The last step is $k_0$th step.

We have $\Theta_0=\tr(a_0,b_0)$, all diagrams $\Theta_k$ lie below $\tr(a_0,b_0)$ and have constant endpoints for $k_0> k>0$. We will argue that the procedure above gives asserted jumps. 

For $\nu(\Theta_0)$ the initial jumps are 1 due to Propositions \ref{propSigma0DoSigmaN} and \ref{propShortEEA}.

For $k>0$ we need only to show that all intermediate polygonal chains are diagrams.

\begin{figure}[h]
\centering
\includegraphics [height=6cm]{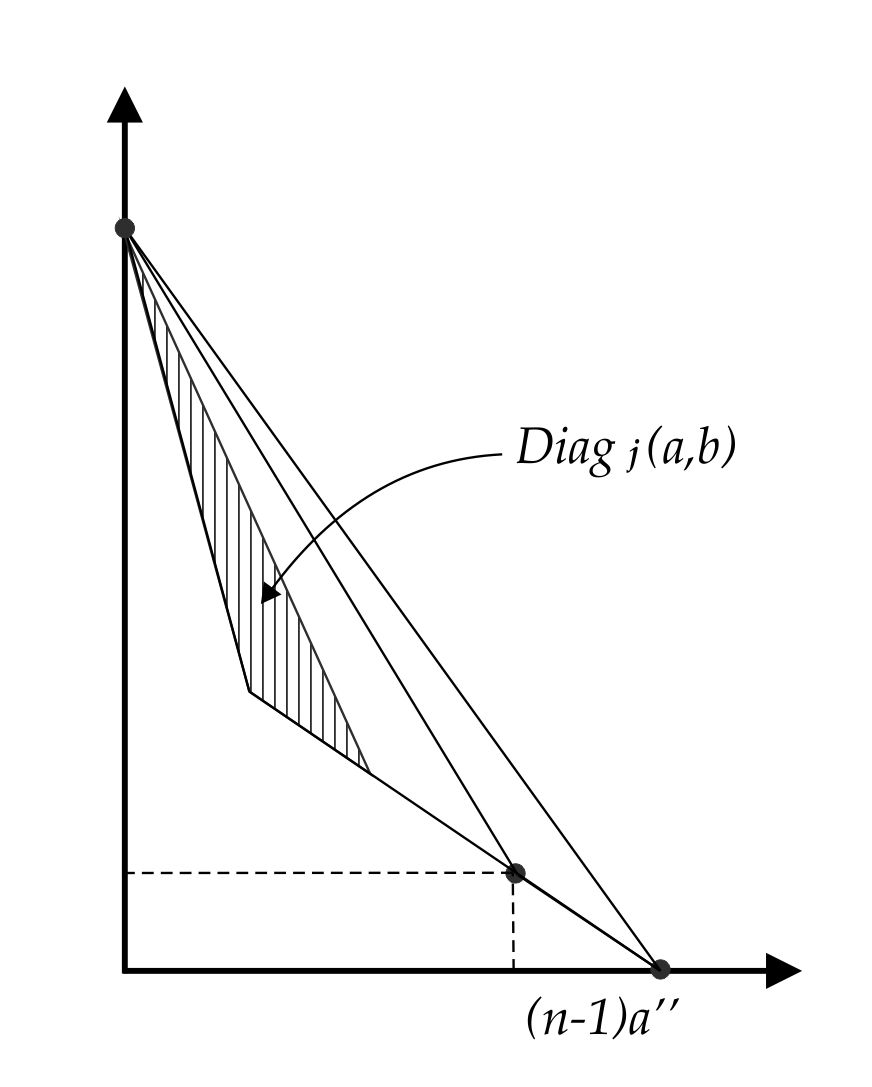} \label{r:3}
\caption{Using Procedure \ref{procka2} for an embedded diagram. Note that $n-1$ is always nonzero and here ${\rm Diag}_j(a,b) = (nj+1)\tr(a',b')+n\tr(a-ja',b-jb')$ and the shaded area is the difference between subsequent diagrams.}
\end{figure}

Indeed, if $a''\neq 0$ recall that due to Remark \ref{rkAllDeformationsAboveGammaN} all points considered for $\tr(p,q)$ lie above $\Gamma^n(\Sigma^{n'-1}) = -{\rm sign}(p,q) ( \left(nn'+1 \right) \tr (a',b') + n\ \tr\left(a'',b''\right) )$. Assume ${\rm sign}(p,q)=-1$. Note that $L\tr(a',b') + \Gamma^n(\Sigma^{n'-1})$ is a diagram with the same endpoints as $\Theta_k$. Hence all intermediate diagrams combined with $L\tr(a',b')$ as the initial segment are diagrams. Recall Property \ref{prtyNewtonNr}. By Lemma \ref{lemTheta} we get that $\nu(\Theta_k)$ has been already attained in the sequence. Hence from Proposition \ref{propSigma0DoSigmaN} the jumps are at most 1. The same argument applies when ${\rm sign}(p,q)=1$.

Moreover, if $a''=0$ or $a'=0$ the same argument gives that from Proposition \ref{propShortEEA} follows that the jumps are at most 1.

Now we only have to compute the total number of jumps. The last diagram due to Remark \ref{rkLastDiagramShort} is $M\tr(1,m+1)+N\tr(1,m)$ for some positive integers $M,N$. Hence we have obtained all numbers ranging from $\nu(\tr(a_0,b_0))$ up to $\nu(M\tr(1,m+1)+N\tr(1,m))$, the difference is double the area and is equal to $r(a_0-r)$. Indeed, we have $a_0=M+N$ and $b_0=M(m+1)+Nm=m(M+N)+M$. Hence $M=r$ and $N=a_0-r$. Thus double the area of the difference is $a_0b_0-Mb_0-a_0mN=N(b_0-a_0m)=r(a_0-r)$. This gives the claim.\hfill$\blacksquare$\\

Note that the above can also be computed explicitly using EEA and inductive definition of $z_k$.

\begin{rk}
Theorem \ref{thmMainCombinatorially} is at its weakest for $q\equiv \pm 1 ({\rm mod}\ p)$, when the function $r(p-r)$ minimises and is equal $p-1$.
\end{rk}

\begin{ex}\label{exDalej}
We continue Example \ref{przyEEA}. For the diagrem $\tr(40,73)$ from Theorem~\ref{thmMainCombinatorially} we get that the sequence of jumps of Newton numbers begins with $33\cdot(40-33)=231$ ones.
\end{ex}

\section{Remarks}\label{sectionRks}

We will indicate one possible use of the algorithm described in this paper in finding all Milnor numbers attained by deformations. Note that this combinatorial approach gives also the form of deformations that have the supposed Milnor number.

Let us look at a continuation of Example~\ref{przyEEA}.
\begin{ex}\label{rk_ex}
Take an irreducible singularity $f$ of the form~\eqref{eqPostacfQSH} with $p=40$ and $q=73$. We claim that all positive integers less than $\mu(f)$ are attained as Milnor numbers of deformations of $f$ i.e. the sequence of jumps is constantly equal $1$.

Indeed, as was already indicated in Example~\ref{exDalej}, we have at least $231$ initial ones in the sequence of jumps of Milnor numbers. Take nondegenerate deformations $F_{k,l}$ of $f$ such that 
$$\Gamma(F_{k,l})=\Gamma(f)+((0,k),(l,0)).$$
Note that the diagram $\Gamma(F_{k,l})$ consists of a single segment. Consider for instance $F_{37,73}$. We have $37$ and $73$ are coprime, moreover $\mu(F_{37,73})>\mu(f)-231>\mu(F_{37,73})-36\cdot(37-36)$. Using Theorem~\ref{mainThm}, we get that the sequence of jumps equal to $1$ is at least as long as $252$. This improves the previous result. 

In the same manner consider deformations $F_{k,l}$ with $(k,l)$ consecutively equal to 
$$
\begin{array}{l}
(39,73),\ (38,73),\ (37,73),\\
(37,73),\ (37,71),\ \dots,\ (37,41)
\end{array}
$$
and apply Theorem~\ref{mainThm} to each. Now one can continue with deformations with $(k,l)$ equal to
$$
\begin{array}{l}
(37,41), \ (36,41),\ \dots,\ (23,41),\\
(23,41),\ (23,40),\ \dots,\ (23,29),\\
(23,29),\ (22,29),\ \dots\ etc
\end{array}
$$
or use the main result of~\cite{BKW} for $k=40$. Precisely, the result we are referring to states that for a homogeneous nondegenerate isolated singularity $f_k$ of degree $k$ all positive integers less than $\mu(f_k)-k+2$ are attained as Milnor numbers of deformations. Note that $\mu(f_{40})-40+2>\mu(F_{37,41})-4\cdot(37-4)$.

Both approaches give the assertion of Example~\ref{rk_ex}. 
\end{ex}

Do note that the computation by hand presented above (which is also easy to implement as a program) is essentially better than straightforward numerical computation. Our numerical experiments with naive algorithms have lasted for hours in the case of the singularity from Example~\ref{rk_ex}, whereas doing it by hand using Theorem~\ref{mainThm} is a matter of minutes. 

Easy generalisation of the above Example is

\begin{cor}\label{cor_ex}
Take an isolated singularity $f$ of the form~\eqref{eqPostacfQSH} with $p<q$ coprime. Suppose there exists an injective sequence of coprime numbers $(p_s,q_s)_{s=1,\dots,v}$ such that $p_s\leq q_s$, both sequences $(p_s), (q_s)$ are non-increasing and
\begin{eqnarray}\label{eqRk_ex}
(p_{s}-1)(q_{s}-1)-r_{s}(p_{s}-r_{s})\le (p_{s+1}-1)(q_{s+1}-1),
\end{eqnarray}
where we denote by $r_s$ the positive integer such that $q_s\equiv r_s ({\rm mod}\ p_s)$.

Then all positive integers between $\mu(f)$ and $(p_{v}-1)(q_{v}-1)-r_{v}(p_{v}-r_{v})$ are attained as Milnor numbers of deformations of $f$. 

Moreover, if 
\begin{equation}\label{eqRk_ex2}
(p_{v}-1)(q_{v}-1)-r_{v}(p_{v}-r_{v})< (p-1)(p-2)+1,
\end{equation}
then all positive integers are attained as Milnor numbers of deformations of~$f$.
\end{cor}

\noindent{\bf Proof.} As in Example~\ref{rk_ex} consider nondegenerate deformations $F_{k,l}$ of $f$ such that 
$$\Gamma(F_{k,l})=\Gamma(f)+((0,k),(l,0)).$$
Since $p_{s},q_{s}$ are coprime, for each deformation $F_{p_{s},q_{s}}$ use Theorem~\ref{mainThm}. A deformation of a deformation is a deformation (can be chosen as a one-parameter deformation as well as from nondegeneracy of deformations one can assume it is nondegenerate) due to the form of the diagrams. Hence we get deformations of $f$ giving Milnor numbers from $(p_{s}-1)(q_{s}-1)$ to $(p_{s}-1)(q_{s}-1)-r_{s}(p_{s}-r_{s})$. The inequality \eqref{eqRk_ex} guaranties that the Milnor number of $F_{p_{s+1},q_{s+1}}$ is greater than $\mu(F_{p_{s},q_{s}})-r_{s}(p_{s}-r_{s})$. Hence all integers between $\mu(f)$ and $(p_{v}-1)(q_{v}-1)-r_{v}(p_{v}-r_{v})$ are attained as Milnor numbers of deformations of $f$.

Moreover, if inequality \eqref{eqRk_ex2} holds, it means that the Milnor number of the deformation $F_{p,p}$ of $f$ is bigger by at least $p-2$ than the last Milnor number already attained. Hence we can use the result that for a homogeneous nondegenerate isolated singularity of degree $p$ all positive integers less or equal $\mu(F_{p,p})-p+2$ are attained as Milnor numbers of deformations from \cite{BKW}. This ends the proof.~\hfill$\blacksquare$

The procedure above may be stated constructively but we cannot at the moment guarantee that we can choose sequences which satisfy inequality~\eqref{eqRk_ex2}. To the contrary, for $p,q$ relatively small, for instance $(5,7)$, such a sequence may not exist (compare \cite{W2}). The question on when such a sequence exists could be possibly resolved using distribution of primes (compare the explicit sequence from Remark~\ref{rk_ex}).

This paper answers in particular to open questions posed in the article \cite{Bo} which stress for constructive methods. Our combinatorial approach in the spirit of previous sections is very powerful in answering these questions. The authors have results also for all bivariate semi-quasi-homogeneous singularities, in particular extending results of this paper on irreducible germs, but we defer the details to a subsequent publication.

\section*{Aknowledgements} Authors were supported by grant NCN 2013/09/D/ST1/03701.

\bibliographystyle{alpha}
\bibliography{bibliografia}

\begin{thebibliography}{BKW14}

\bibitem[Arn04]{ArnoldProblems}
Vladimir~I. Arnold.
\newblock {\em Arnold's problems}.
\newblock Springer-Verlag, Berlin; PHASIS, Moscow, 2004.
\newblock Translated and revised edition of the 2000 Russian original, With a
  preface by V. Philippov, A. Yakivchik and M. Peters.

\bibitem[BK14]{BK}
Szymon Brzostowski and Tadeusz Krasi{\'n}ski.
\newblock The jump of the {M}ilnor number in the {$X_9$} singularity class.
\newblock {\em Cent. Eur. J. Math.}, 12(3):429--435, 2014.

\bibitem[BKW14]{BKW}
S.~{Brzostowski}, T.~{Krasinski}, and J.~{Walewska}.
\newblock {Milnor numbers in deformations of homogeneous singularities}.
\newblock {\em ArXiv}, 2014.

\bibitem[Bod07]{Bo}
Arnaud Bodin.
\newblock Jump of {M}ilnor numbers.
\newblock {\em Bull. Braz. Math. Soc. (N.S.)}, 38(3):389--396, 2007.

\bibitem[GLS07]{GreuelShustinL}
G.-M. Greuel, C.~Lossen, and E.~Shustin.
\newblock {\em Introduction to singularities and deformations}.
\newblock Springer Monographs in Mathematics. Springer, Berlin, 2007.

\bibitem[GN12]{GreuelNonDegN2}
Gert-Martin Greuel and Hong~Duc Nguyen.
\newblock Some remarks on the planar {K}ouchnirenko's theorem.
\newblock {\em Rev. Mat. Complut.}, 25(2):557--579, 2012.

\bibitem[GZ93]{GuseinZade}
S.~M. Guse{\u\i}n-Zade.
\newblock On singularities that admit splitting off {$A_1$}.
\newblock {\em Funktsional. Anal. i Prilozhen.}, 27(1):68--71, 1993.

\bibitem[Kou76]{Kush}
A.~G. Kouchnirenko.
\newblock Poly\`edres de {N}ewton et nombres de {M}ilnor.
\newblock {\em Invent. Math.}, 32(1):1--31, 1976.

\bibitem[Len08]{Len08}
Andrzej Lenarcik.
\newblock On the {J}acobian {N}ewton polygon of plane curve singularities.
\newblock {\em Manuscripta Math.}, 125(3):309--324, 2008.

\bibitem[P{\l}o14]{Pl14}
Arkadiusz P{\l}oski.
\newblock A bound for the {M}ilnor number of plane curve singularities.
\newblock {\em Cent. Eur. J. Math.}, 12(5):688--693, 2014.

\bibitem[Wal99]{WallNonDeg}
C.~T.~C. Wall.
\newblock Newton polytopes and non-degeneracy.
\newblock {\em J. Reine Angew. Math.}, 509:1--19, 1999.

\bibitem[Wal10]{W1}
Justyna Walewska.
\newblock The second jump of {M}ilnor numbers.
\newblock {\em Demonstratio Math.}, 43(2):361--374, 2010.

\bibitem[Wal13]{W2}
Justyna Walewska.
\newblock Jumps of the {M}ilnor numbers in families of non-degenerate and
  non-convenient singularities.
\newblock In {\em Analytic and Algebraic Geometry}, Proceedings of Conference
  on Analytic and Algebraic Geometry. 2013.

\end{thebibliography}

\end{document}